\magnification=1100
\input amssym.def
\input amssym.tex
\topskip=1cm

\advance\hoffset by 2 truecm
\vsize= 22,5 truecm
\hsize= 14 truecm
\def\makefootline{\baselineskip=72pt\line{\the\footline}} \overfullrule=0pt



\centerline{\bf On Linearity of Nonclassical Differentiation } \bigskip\medskip

\centerline{\bf Serguei SAMBORSKI}
\smallskip
\centerline{Universit\'{e} de Caen, Math\'{e}matiques, 14032 CAEN Cedex, France}

\medskip
\centerline{\vbox{\hrule width 1,5cm}}

\vskip1,5cm
\parindent=0,5cm

\noindent{\bf Abstract}
\smallskip We introduce a real vector space composed of set-valued maps on an open set X and note it by S. It is a complete
metric space and a conditionally complete lattice. The set of continuous functions on X is dense in S as in a metric space
and as in a lattice. Thus the constructed space plays the same role for the space of continuous functions with uniform
convergence as the field of reals plays for the field of rationals. The classical gradient may be extended in the space S as
a closed operator. If a function f belongs to the domain of this extension, then f is locally lipschitzian and the values
of our gradient coincide with the values of Clarke's gradient. However, unlike Clarke's gradient, our generalized gradient
is a linear operator. \bigskip
\noindent{\bf Key words}\smallskip functional metric spaces, functional lattices, extensions of differentiation, Clarke's
gradient, quasi continuous functions. \bigskip \noindent{\bf AMS Subject Classification}\smallskip 26A24, 28A15, 46E05, 54C35
\bigskip
\noindent{\bf Author}\smallskip\ Dr. Prof. Samborski Serguei\par Mathematiques,~Universite de CAEN,~14032 CAEN
Cedex,~France\par e-mail: samborsk@math.unicaen.fr\par fax : +33~2~31~56~73~20
\bigskip\bigskip\break
\noindent{\bf Sorry, this article is being rewritten. Please email the author to be informed about its availability.}

\end